\documentclass{article}

\usepackage{proof}
\usepackage{latexsym}
\usepackage{amsfonts}
\usepackage{amssymb}
\usepackage{cite}

\newcommand{\Natural}{\mathbb{N}}

\newcommand{\brem}{\begin{remark}}
\newcommand{\erem}{\end{remark}}

\newcommand{\blem}{\begin{lemma}}
\newcommand{\elem}{\end{lemma}}
\newcommand{\bth}{\begin{theorem}}
\newcommand{\ethm}{\end{theorem}}
\newcommand{\benu}{\begin{enumerate}}
\newcommand{\eenu}{\end{enumerate}}
\newcommand{\bdes}{\begin{description}}
\newcommand{\edes}{\end{description}}
\newcommand{\bdf}{\begin{definition}}
\newcommand{\edf}{\end{definition}}
\newcommand{\bcor}{\begin{cor}}
\newcommand{\ecor}{\end{cor}}
\newcommand{\bprp}{\begin{proposition}}
\newcommand{\eprp}{\end{proposition}}
\newcommand{\bmlem}{\begin{mlemma}}
\newcommand{\emlem}{\end{mlemma}}
\newcommand{\bclm}{\begin{claim}}
\newcommand{\eclm}{\end{claim}}
\newcommand{\bprf}{{\bf Proof}.\hspace{2mm}}

\newcommand{\eprf}{\hspace*{\fill} $\Box$}

\newcommand{\beqn}{\begin{equation}}
\newcommand{\eeqn}{\end{equation}}
\newcommand{\beqnarr}{\begin{eqnarray}}
\newcommand{\eeqnarr}{\end{eqnarray}}
\newcommand{\beqnarrs}{\begin{eqnarray*}}
\newcommand{\eeqnarrs}{\end{eqnarray*}}

\newcommand{\spand}{\,\&\,}

\newtheorem{theorem}{Theorem}[section]
\newtheorem{definition}[theorem]{Definition}
\newtheorem{proposition}[theorem]{Proposition}
\newtheorem{lemma}[theorem]{Lemma}
\newtheorem{cor}[theorem]{Corollary}
\newtheorem{mlemma}[theorem]{Main Lemma}
\newtheorem{claim}[theorem]{Claim}

\newtheorem{remark}[theorem]{Remark}

\newcommand{\alp}{\alpha}

\newcommand{\veps}{\varepsilon}
\newcommand{\del}{\delta}
\newcommand{\Del}{\Delta}
\newcommand{\ome}{\omega}

\newcommand{\bet}{\beta}
\newcommand{\gam}{\gamma}
\newcommand{\Gam}{\Gamma}

\newcommand{\sig}{\sigma}
\newcommand{\Sig}{\Sigma}
\newcommand{\tht}{\theta}

\newcommand{\lam}{\lambda}

\newcommand{\vphi}{\varphi}

\newcommand{\fal}{\forall}
\newcommand{\exi}{\exists}

\newcommand{\Rarw }{\Rightarrow}

\newcommand{\lrarw}{\leftrightarrow}
\newcommand{\Lrarw}{\Leftrightarrow}

\newcommand{\la}{\langle}
\newcommand{\ra}{\rangle}

\title{
Two remarks on proof theory of first-order arithmetic
}
\author{Toshiyasu Arai 
\\
Graduate School of Mathematical Sciences,
University of Tokyo,
\\
3-8-1 Komaba, Meguro-ku,
Tokyo 153-8914, JAPAN
\\
tosarai@ms.u-tokyo.ac.jp}
\date{}
\begin{document}
\maketitle

\begin{abstract}
In this note let us give two remarks on
proof-theory of {\sf PA}.
First a derivability relation is introduced to bound witnesses for provable $\Sig_{1}$-formulas in {\sf PA}.
Second Paris-Harrington's proof for their independence result is reformulated
to show a `consistency' proof of {\sf PA} based on a combinatorial principle.
\end{abstract}

\section{A derivability relation to bound witnesses}

\begin{theorem} \label{th:BW}
Each provably computable function in {\sf PA}
is dominated by a Hardy function $H_{\alp}\,(\alp<\veps_{0})$.
\end{theorem}
Theorem \ref{th:BW} is a classic result in proof theory of {\sf PA}, and
there are known 
several proofs of it by G. Kreisel, S. Wainer, H. Schwichtenberg, et al.
Let us introduce a derivability relation to bound witnesses for provable $\Sig_{1}$-formulas in {\sf PA},
which is a variant of one given in Chapter 4 of Schwichtenberg-Wainer\cite{SW}.
The latter seems to be inspired from the proof in
Buchholz-Wainer\cite{BW}.

In this section by ordinals we mean ordinals$<\veps_{0}$.
When $\gam+\alp=\gam\#\alp$, we write $\gam\dot{+}\alp$ for $\gam+\alp$.
In this case we say that $\gam\dot{+}\alp$ is defined.
Also let $0\dot{+}\gam=\gam\dot{+}0:=\gam$.

Fundamental sequences $\{\lam[x]\}_{x\in\Natural}$ for limit ordinals $\lam<\veps_{0}$ are defined as 
follows.
Let
$0[x]=0$, $(\alp+1)[x]=\alp$.
For ordinals $0<\alp<\veps_{0}$ in Cantor normal form
$\alp=\ome^{\alp_{0}}m_{0}+\cdots+\ome^{\alp_{k}}m_{k}
\,(\alp_{0}>\cdots>\alp_{k}\geq 0,0<m_{0},\ldots,m_{k}<\ome)$,
let
\[
\alp[x]=\left\{
\begin{array}{ll}
\ome^{\alp_{0}}m_{0}+\cdots+\ome^{\alp_{k}}(m_{k}-1)+\ome^{\alp_{k}[x]}  
& \alp_{k}\mbox{ is a limit ordinal}
\\
\ome^{\alp_{0}}m_{0}+\cdots+\ome^{\alp_{k}}(m_{k}-1)+\ome^{\bet}x
& \alp_{k}=\bet+1 
\\
\ome^{\alp_{0}}m_{0}+\cdots+\ome^{\alp_{k}}(m_{k}-1)
& \alp_{k}=0 
\end{array}
\right.
\]

$\alp<_{n}\bet$ denotes the transitive closure of the relation
$\{(\alp,\bet): \alp=\bet[n]\}$.
Let $\alp\leq_{n}\bet:\Lrarw \alp<_{n}\bet\lor\alp=\bet$.

Hardy functions $H_{\alp}$ are defined by transfinite recursion on ordinals $\alp<\veps_{0}$:
$H_{0}(x)=x$, $H_{\alp+1}(x)=H_{\alp}(x+1)$, and
$H_{\lam}(x)=H_{\lam[x]}(x)$ for limit ordinals $\lam$.

{\sf PA} denotes the first-order arithmetic in the language $\{0,1,+,\cdot, \lam x,y.\,x^{y}, <,=\}$.
${\rm dg}(A)<\ome$ for formulas $A$ is defined by
${\rm dg}(A)=0$ if $A\in\Del_{0}$, 
${\rm dg}(A_{0}\circ A_{1})=\max\{{\rm dg}(A_{i}):i<2\}\,(\circ\in\{\lor,\land\})$,
and ${\rm dg}(Qx A)={\rm dg}(A)+1\, (Q\in\{\exi,\fal\})$.

\bdf
{\rm
For ordinals $\gam,\alp<\veps_{0}$, natural numbers $k,c<\ome$ and 
finite sets $\Gam$ of sentences,
a derivability relation
$(\gam, k)\vdash^{\alp}_{c}\Gam$ is defined by recursion on $\alp$ as follows.
Let $\ell=H_{\gam}(k)$.

\begin{description}
\item [(Ax)] 
If $\Gamma$ contains a true $\Del_{0}$-sentence, denoted by $\top$ ambiguously,
then
$(\gam, k)\vdash^{\alp}_{c}\Gam$ holds for any $k,c<\ome$ and $\gam, \alp<\veps_{0}$.
 
\item [(cut)] 
If
$\alp_{0}+1\leq_{\ell}\alp$,
$(\gam, k)\vdash^{\alpha_{0}}_{c}\Gamma,A$ with ${\rm dg}(A)<c$, and
$(\gam, k)\vdash^{\alpha_{0}}_{c} \Gamma,\lnot A$, then
$(\gam, k)\vdash^{\alp}_{c}\Gam$.


\item [($\land$)]  
If
$\alp_{0}<_{\ell}\alp$, $(A_{0}\land A_{1})\not\in\Del_{0}$,
$(A_{0}\land A_{1}) \in \Gamma$, and
$(\gam, k)\vdash^{\alpha_{0}}_{c}\Gamma, A_{i}$ for any $i=0,1$,
then,
$(\gam, k)\vdash^{\alp}_{c}\Gam$.

\item [($\lor$)] 
If
$\alp_{0}+1\leq_{\ell}\alp$, $(A_{0}\lor A_{1})\not\in\Del_{0}$,
$(A_{0}\lor A_{1}) \in \Gamma$, and
$(\gam, k)\vdash^{\alpha_{0}}_{c}\Gamma, A_{i}$ for an $i=0,1$, then
$(\gam, k)\vdash^{\alp}_{c}\Gam$.

\item [($\forall\ome$)]  
If
$\alp_{0}<_{\ell}\alp$, $\forall x A(x)\not\in\Del_{0}$,
$\forall x A(x) \in \Gamma$, and 
$(\gam, \max\{k,n\})\vdash^{\alpha}_{c}\Gamma, A(n)$ for any $n\in\Natural$, then
$(\gam, k)\vdash^{\alp}_{c}\Gam$.

\item [($\exists$)] 
If
$\alp_{0}+1\leq_{\ell}\alp$, $\exi x A(x)\not\in\Del_{0}$,
$\exists x A(x) \in \Gamma$, and there exists a natural number $n$ such that
$(\gam, k)\vdash^{\alpha_{0}}_{c} A(n),\Gamma$ and
$n< \ell$, then
$(\gam, k)\vdash^{\alp}_{c}\Gam$.

\end{description}
}
\edf

For natural numbers $k$ and $\Sig_{1}$-sentences  $\exi x A\, (A\in\Del_{0})$,
$k\models\exi x A:\Lrarw \Natural\models\exi x<k\, A$.
For finite sets $\Gam$ of $\Sig_{1}$-sentences,
$k\models\Gam$ iff there exists a $B\in\Gam$ such that $k\models B$.

\blem\label{lem:BWbounding}
{\rm (Bounding)}
Let $\Gam$ be a  finite set of $\Sig_{1}$-sentences.
\\
Assume $(\gam, k)\vdash^{\alp}_{0}\Gam$.
Then $H_{\gam}(k)\models\Gam$.
\elem
\bprf
This is seen by induction on ordinals $\alp$.
Consider the case when $(\gam, k)\vdash^{\alpha}_{0} \Gamma$ follows
from an inference rule $(\exi)$.
Let $\exi x A(x)$ be its principal (major) formula.
We have $\alpha_{0}+1\leq_{k}\alpha$ and
$(\gam, k)\vdash^{\alpha_{0}}_{0} A(m),\Gamma$ with $m< H_{\gam}(k)$.
IH yields $H_{\gam}(k)\models\{A(m)\}\cup\Gam$.
\eprf

\begin{lemma}\label{lem:BWinv}
\begin{enumerate}
\item\label{lem:BWinv1}
 {\rm (Weakening)}
Let $(\gam, k)\vdash^{\alpha}_{c}\Gamma$, and assume
$\fal n[H_{\gam}(\max\{k,n\})\leq H_{\gam_{1}}(\max\{k_{1},n\})]$.
Then for $\Gamma\subset\Gamma_{1}$, 
$\alp\leq_{H_{\gam_{1}}(k_{1})}\alp_{1}$, and $c \leq c_{1}$, 
$(\gam_{1}, k_{1})\vdash^{\alpha_{1}} _{c_{1}}\Gamma_{1}$ holds.

\item\label{lem:BWinv2}
{\rm (False)}
For a false $\Del_{0}$-sentence $\bot$, if
$(\gam, k)\vdash^{\alpha}_{c}\Gamma,\bot$, then
$(\gam, k)\vdash^{\alpha}_{c}\Gamma$.

\item\label{lem:BWinv3}
 {\rm (Inversion)}
 If
$(\gam, k)\vdash^{\alpha}_{c}\Gamma,\forall xA(x)$, then
$(\gam, \max\{k,n\})\vdash^{\alpha}_{c} \Gamma,A(n)$ for any $n\in\Natural$.

\item\label{lem:BWinv4}
Assume $\gam\dot{+}\del$ is defined,
$n<H_{\del}(k)$ and $(\gam,\max\{k,n\})\vdash^{\alp}_{c}\Gam$.

Then
$(\gam\dot{+}\del,k)\vdash^{\alp}_{c}\Gam$.
\eenu
\end{lemma}
\bprf
\ref{lem:BWinv}.\ref{lem:BWinv2}.
No $\Del_{0}$-formula is a principal formula of an inference rule.
\\
\ref{lem:BWinv}.\ref{lem:BWinv4}.
This is seen from Lemma \ref{lem:BWinv}.\ref{lem:BWinv1}.
\eprf

\begin{lemma}\label{lem:BWreduction}
{\rm (Reduction)}
Let $C$ be a sentence such that ${\rm dg}(C)\leq c$ and
$C$ is either one of the form $\fal x\,A, A_{0}\land A_{1}\not\in\Del_{0}$
or a true $\Del_{0}$-sentence $\top$.
Assume that $\alp\dot{+}\bet$ is defined, and $\gam=\ome^{\gam_{0}}m$.

Suppose $(\gam, k)\vdash^{\alp}_{c}\Gam,C$
and $(\gam,k)\vdash^{\bet}_{c}\lnot C,\Del$.
Then
$(\gam\cdot 2,k)\vdash^{\alp\dot{+}\bet}_{c}\Gam,\Del$.
\end{lemma}
\bprf
By induction on $\bet$.
Consider the case when $(\gam,k)\vdash^{\bet}_{c}\lnot C,\Del$ follows from an
inference rule $(\exi)$ with its principal formula $\lnot C\equiv(\exi x\lnot A(x))$.
There are $\bet_{0}+1\leq_{\ell}\bet$, $n< \ell=H_{\gam}(k)$ for which
$(\gam, k)\vdash^{\bet_{0}}_{c}\Del,\lnot C, \lnot A(n)$.
IH yields
$(\gam\cdot 2,k)\vdash^{\alp\dot{+}\bet_{0}}_{c}\Gam,\Del, \lnot A(n)$.
On the other hand we have 
$(\gam,\max\{k,n\})\vdash^{\alp}_{c}\Gam,A(n)$ by Lemma \ref{lem:BWinv}.\ref{lem:BWinv3}.
Lemma \ref{lem:BWinv}.\ref{lem:BWinv4} with $n< \ell$ yields
$(\gam\cdot 2k)\vdash^{\alp}_{c}\Gam,A(n)$.
Since ${\rm dg}(A(n))<c$, we obtain
$(\gam\cdot 2,k)\vdash^{\alp\dot{+}\bet}_{c}\Gam,\Del$ by a $(cut)$.
\eprf

\begin{lemma}\label{lem:BWelim}
{\rm (Elimination)}
Assume that $\gam\dot{+}\alp$ is defined for $\gam=\ome^{\gam_{0}}m\geq\ome$, and $k\geq 2$.
If $(\gam,k)\vdash^{\alpha}_{c+1} \Gamma$, then
$(\ome^{\gam\dot{+}\alp}+\gam,k)\vdash^{\omega^{\alpha}}_{c} \Gamma$.
\end{lemma}
\bprf
By induction on $\alp$.
Suppose that $(\gam,k)\vdash^{\alpha}_{c+1} \Gamma$ follows from an inference rule $I$.
Let $\ell=H_{\gam}(k)$.

First consider the case when $I$ is an $(\exi)$.
We have $(\gam,k)\vdash^{\bet}_{c+1}\Gam,B(n)$ for
$\bet+1\leq_{\ell}\alp$, $n< \ell$, and $(\exi x B(x))\in\Gam$.
IH yields
$(\ome^{\gam\dot{+}\bet}+\gam,k)\vdash^{\ome^{\bet}}_{c}\Gam,B(n)$.
By $\bet+1\leq_{\ell}\alp$,
we obtain  
$\ome^{\bet}+1<_{k^{\prime}}\ome^{\alp}$ and
$\ome^{\gam\dot{+}\bet}<_{k^{\prime}}\ome^{\gam\dot{+}\alp}$
for any $k^{\prime}\geq H_{\gam}(k)\geq 2$.
Hence for $n\geq k$ we obtain
$H_{\ome^{\gam\dot{+}\bet}}(H_{\gam}(n))\leq H_{\ome^{\gam\dot{+}\alp}}(H_{\gam}(n))$,
and
$(\ome^{\gam\dot{+}\alp}+\gam,k)\vdash^{\ome^{\bet}}_{c}\Gam,B(n)$ 
by Lemma \ref{lem:BWinv}.\ref{lem:BWinv1}.
An $(\exi)$ yields $(\ome^{\gam\dot{+}\alp}+\gam,k)\vdash^{\ome^{\alp}}_{c}\Gam$.

Next consider the case when $I$ is a $(\fal\ome)$.
We have $(\gam,\max\{k,n\})\vdash^{\bet}_{c+1}\Gam,B(n)$ for
$\bet<_{\ell}\alp$ and $(\fal x B(x))\in\Gam$.
IH yields
$(\ome^{\gam\dot{+}\bet}+\gam,\max\{k,n\})\vdash^{\ome^{\bet}}_{c}\Gam,B(n)$.
We see
$(\ome^{\gam\dot{+}\alp}+\gam,\max\{k,n\})\vdash^{\ome^{\bet}}_{c}\Gam,B(n)$
similarly for the case $(\exi)$.
A $(\fal\ome)$ yields $(\ome^{\gam\dot{+}\alp}+\gam,k)\vdash^{\ome^{\alp}}_{c}\Gam$.

Finally consider the case when $I$ is a $(cut)$ with the cut formula $C$.
We have ${\rm dg}(C)\leq c$,
$(\gam,k)\vdash^{\bet}_{c+1}\Gam,C$ and
$(\gam,k)\vdash^{\bet}_{c+1}\Gam,\lnot C$ for $\bet+1\leq_{\ell}\alp$.
IH yields
$(\ome^{\gam\dot{+}\bet}+\gam,k)\vdash^{\ome^{\bet}}_{c}\Gam,C$ and
$(\ome^{\gam\dot{+}\bet}+\gam,k)\vdash^{\ome^{\bet}}_{c}\Gam,\lnot C$.
We have
$H_{\gam}(k^{\prime})<H_{\ome^{\gam\dot{+}\bet}}(k^{\prime})$ for $k^{\prime}\geq 2$.
Lemma \ref{lem:BWinv}.\ref{lem:BWinv1} yields for $k\geq 2$,
$(\ome^{\gam\dot{+}\bet}\cdot 2,k) \vdash^{\ome^{\bet}}_{c}\Gam,C$ and
$(\ome^{\gam\dot{+}\bet}\cdot 2,k)\vdash^{\ome^{\bet}}_{c}\Gam,\lnot C$.
Lemma  \ref{lem:BWreduction} yields
$(\ome^{\gam\dot{+}\bet}\cdot 4,k)\vdash^{\ome^{\bet}\cdot 2}_{c}\Gam$.
On the other hand we have
$\ome^{\gam\dot{+}\bet}\cdot 4<_{H_{\gam}(k^{\prime})}\ome^{\gam\dot{+}\alp}$
and
$\ome^{\bet}\cdot 2<_{H_{\gam}(k)}\ome^{\alp}$
for $k^{\prime}\geq k$ with
$H_{\gam}(k)\geq H_{\ome}(k)=2k\geq 4$, $k\geq 2=|\ome|$, and $\gam\geq\ome$.
This yields
$H_{\ome^{\gam\dot{+}\bet}\cdot 4}(k^{\prime})\leq
H_{\ome^{\gam\dot{+}\bet}\cdot 4}(H_{\gam}(k^{\prime}))
<
H_{\ome^{\gam\dot{+}\alp}}(H_{\gam}(k^{\prime}))$.
Lemma \ref{lem:BWinv}.\ref{lem:BWinv1} yields
$(\ome^{\gam\dot{+}\alp}+\gam,k)\vdash^{\ome^{\alp}}_{c}\Gam$.
\eprf

\blem\label{lem:BWembed}{\rm (Embedding)}
Let $\Gam(\vec{a})$ be a sequent with free variables $\vec{a}$.
$\Gam(\vec{n})$ dnotes the result of replacing variables $\vec{a}$
by numerals $\vec{n}$.
If ${\sf PA}\vdash\Gam(\vec{a})$, then there are natural numbers
$\ell,d,m,c<\ome$ such that
$(\ome^{2}\ell, \max(\{1\}\cup\vec{n}))\vdash^{\ome\cdot d+m}_{c}\Gam(\vec{n})$
for any $\vec{n}$.
\elem
\bprf
Let $n^{\prime}=\max(\{1\}\cup\vec{n})$.
Consider the case when $\Gam$ follows from an $(\exi)$ with its principal formula
$(\exi x\,A(x))\in\Gam$.
\[
\infer[(\exi)]{\Gam}{\Gam, A(t)}
\]
IH yields $\ome^{2}\ell,n^{\prime}\vdash^{\ome\cdot d+m}_{c}\Gam, A(t(\vec{n}))$.
Pick a $k$ so that $t(\vec{n})<H_{\ome^{2}k}(n^{\prime})$ with $H_{\ome^{2}}(x)\geq 2^{x}$.
An inference rule $(\exi)$ with $\ell_{1}=\max\{\ell,k\}$ yields
$(\ome^{2}\ell_{1},n^{\prime})\vdash^{\ome\cdot d+m+1}_{c}\Gam$.
\eprf
\\

Let us prove Theorem \ref{th:BW}.
Assume 
${\sf PA}\vdash \exists y\, \tht(x,y)\,(\tht\in\Del_{0})$.
By Lemma \ref{lem:BWembed} pick natural numbers $\ell>0,d,m,c$ so that for any
$n\in\Natural$,
$\ome^{2}\ell,\max\{1,n\}\vdash^{\ome\cdot d +m}_{c} \exi y\, A(n,y)$.
Define ordinals $\gam_{c}$ by $\gam_{0}=\ome^{2}\ell$,
$\gam_{c+1}=\ome^{\gam_{c}\dot{+}\omega_{c}(\ome\cdot d+m)}\cdot 2$.
Lemma \ref{lem:BWelim} yields
$\gam_{c},n\vdash^{\alpha}_{0} \exi u\exists y,z A(n,y,z)$
for $n\geq 2$ and $\alpha=\omega_{c}(\ome\cdot d+m)$.
We obtain $\exi y<H_{\gam_{c}}(n)\,A(n,y)$ by Lemma \ref{lem:BWbounding}
for $n\geq 2$.

\brem
{\rm
Our proof does not give the optimal bound $2_{k}(\ome^{2})$
for fragments I$\Sig_{k}\,(k>0)$.
The number $c$ in Lemma \ref{lem:BWembed} is bounded by $k+1$.
Then $\gam_{k+1}<\ome_{k+1}(\ome^{3})$.
To obtain a better bound, the derivability relation $(\gam,k)\vdash^{\alp}_{c}\Gam$
is modified as follows.

First let us assume that every formula is in prenex normal form $\bigcup_{c}(\Sig_{c}\cup\Pi_{c})$, and 
${\rm dg}_{1}(A)= \min\{c:A\in\Sigma_{c+1}\cup\Pi_{c+1}\}$.
$(\gam,k)\vdash^{\alpha}_{c} \Gamma$ is defined through ${\rm dg}_{1}(A)$.

Elimination Lemma \ref{lem:BWelim} is stated as follows:
Assume $(\gam,k)\vdash^{\alpha}_{c+1} \Gamma$ 
for $k,c\geq 1$.
Then $(\gam\cdot 4^{\alp}+\gam,k)\vdash^{2^{\alpha}}_{c} \Gamma$.

Likewise Bounding Lemma \ref{lem:BWbounding} is stated as follows:
Let $\Gam\subset\Sig_{1}$ and $(\gam,k)\vdash^{\alp}_{1}\Gam$.
Then $H_{\gam\cdot 4^{\alp}}(H_{\gam}(k))\models\Gam$ for $k\geq 1$.

Now assume ${\rm I}\Sig_{k}\vdash \exists y A(x,y)$ for a $\exists y A(x,y) \in \Sigma_{1}$
and $k>0$. 
We see $\fal n\geq 1\exi y<H_{\del_{k}}(n)\, A(n,y)$
for an ordinal $\del_{k}<2_{k}(\ome^{2})$.
}
\erem

\section{A consistency proof of {\sf PA} based on a combinatorial principle}

\bdf
{\rm
\benu
\item
A subset $H\subset X$ is \textit{diagonal homogeneous} for a partition 
\\
$P:[X]^{1+n}\to c$ if
\[
\fal x_{0}\in H\fal a<x_{0}\fal Y,Z\in[H]^{n}[x_{0}<Y \spand x_{0}<Z \Rarw P(\{a\}\cup Y)=P(\{a\}\cup Z)]
\]
\item
For positive integers $n,m,c$, $X\to_{\Del}(m)^{1+n}_{c}$ designates that
for any partition $P:[X]^{1+n}\to c$, there exists a diagonal homogeneous set $H\in[X]^{m}$.
\item
\textit{Diagonal Homogeneous  principle} denoted by {\rm DH} states that
\\
$\fal n,m,c>0\exi K>1+n[K\to_{\Del} (m)^{1+n}_{c}]$.
\item
Let $\Gam=\{\vphi[y,x_{1},\ldots,x_{n}],\ldots\}$ be a set of formulas in variables 
$y,x_{1},\ldots,x_{n}$, and $D\subset K=\{0,1,\ldots,K-1\}$.
$D$ is said to be a \textit{diagonal indiscernibles} with respect to $\Gam$ (and $K$) if
for any
\[
a<i_{0}<
\begin{array}{lllll}
i_{1} & < & \cdots & < & i_{n}
\\
j_{1}& < & \cdots & < & j_{n} 
\end{array}
,\, (a\in K, \, i_{0},i_{1},\ldots,i_{n},j_{1},\ldots,j_{n}\in D)
\]
$\Natural\models \vphi[a,i_{1},\ldots,i_{n}] \lrarw \vphi[a,j_{1},\ldots,j_{n}]$ holds.
\eenu
}
\edf
It is easy to see that the infinite Ramsey theorem together with K\"onig's lemma implies DH.
\bprp\label{th:diagonalindiscernible}
Let 
$\Gam=\{\vphi_{r}[x_{0},x_{1},\ldots,x_{n}] : r<m\}$ be a finite set of formulas in the language of {\sf PA},
and $k$ an integer.

{\rm DH} $(\exi K[K\to_{\Del}(k)^{1+n}_{2^{m}}])$ yields a diagonal indiscernible set
$D=\{a_{1}<\cdots<a_{k}\}\subset K$ with respect to $\Gam$.
\eprp
\bprf
Let
$P:[K]^{1+n}\to 2^{m}$ be the partion
$P(b_{0},b_{1},\ldots,b_{n})=\{i<m: \vphi_{i}[b_{0},b_{1},\ldots,b_{n}]\}$,
and 
$D=\{a_{1}<\cdots<a_{k}\}\subset K$ be a diagonal homogeneous set for the partition $P$.
\eprf
\\

Pari-Harrington's principle PH states that
$\fal n,m,c\exi K>n[K\to_{*}(m)^{n}_{c}]$, where $K\to_{*}(m)^{n}_{c}$
designates that for any partition $P:[K]^{n}\to c$ there exists a homogeneous set $H\subset K$
with $\# H\geq \min H$.
The proof in Paris-Harrington\cite{PH} of the independence of PH from {\sf PA}
consists of two steps.
First ${\rm Con}(T)\to{\rm Con}({\sf PA})$ for an extension $T$ of {\sf PA},
and ${\rm PH}\to{\rm Con}(T)$.
$T$ is obtained from {\sf PA} by adding an infinite list $\{c_{i}\}_{i<\ome}$
of (individual constants intended to denote) diagonal indiscernibles $c_{i}$.
The purely model-theoretic proof of the independence is given in Kanamori-McAloon\cite{KM}.
In these proofs the principle DH is implicit, and crucial.

Let us reformulate the proof of ${\rm Con}(T)\to{\rm Con}({\sf PA})$ in \cite{PH}
to get a purely proof-theoretic result for {\sf PA}, cf.\,Theorem \ref{lem:DHCon}.

\bth\label{th:DHCon}
${\sf EA}+\fal x\exi y(2_{x}=y) \vdash{\rm DH}\to \mbox{{\rm 1-Con}}({\sf PA})$.
\end{theorem}

A formula $\vphi\equiv\left(Q_{1}z_{1}Q_{2}z_{2}\cdots Q_{n}z_{n}\tht\right)$
with a $\Del_{0}$-matrix $\tht$ and alternating quantifiers 
$Q_{1}, Q_{2},\ldots,Q_{n}$ is a $\Sig_{n}$-formula [$\Pi_{n}$-formula] 
if $Q_{1}\equiv \exi$ [$Q_{1}\equiv\fal$], resp.

In this section {\sf PA} is formulated in an applied one-sided sequent calculus.
Besides usual inference rules for first-order logic,
there are two inference rules for complete induction and axioms
for constants.
The inference rule for complete induction is stated as follows.
\beqn
\label{eq:VJrule}
\infer{\Gam}
{\Gam,A(0,\vec{s}) & \Gam,\lnot A(a,\vec{s}), A(a+1,\vec{s}) & \Gam,\lnot A(t,\vec{s})}
\eeqn
where $A\in\bigcup_{n}(\Sig_{n}\cup\Pi_{n})$.

Let $\fal\vec{x}\tht(\vec{x})$ be a $\Pi_{1}$-axiom for constants $0,1,+,\cdot,\lam x,y.\,x^{y},<,=$.
Then for each list $\vec{t}$ of terms, the following is an inference rule in {\sf PA}:
\[
\infer{\Gam}{\lnot \tht(\vec{t}),\Gam}
\]
The applied calculus admits a cut-elimination.
${\sf PA}\vdash^{d}\Gam$ designates that there exists a derivation of the sequent $\Gam$ 
in depth$\leq d$.

In what follows argue in ${\sf EA}+\fal x\exi y(2_{x}=y)$.
Let $\pi$ be a (cut-free) derivation of a $\Sig_{1}$-sentence
$\exi x\,\tht_{0}(x)$ in {\sf PA}.
Each formula in $\pi$ is in $\bigcup_{n}(\Sig_{n}\cup\Pi_{n})$.

For a formula 
$\vphi\equiv\left(Q_{1}z_{1}Q_{2}z_{2}\cdots Q_{n}z_{n}\tht\right)\in\Sig_{n}\cup\Pi_{n}$, 
let $q(\vphi)=n$.
$q(\tht)=0$ for $\tht\in\Del_{0}$.
Let $Fml(\pi)$ denote the set of all formulas $\vphi$
appearing in $\pi$.
Then $q(\pi):=\max\{q(\vphi): \vphi\in Fml(\pi)\}$.

Second
$m(\pi)=\#\left(\{\vphi\in Fml(\pi) : q(\vphi)>0\}\right)$.
Third $d(\pi)$ denotes the depth of $\pi$:
${\sf PA}\vdash^{d(\pi)}\exi x\,\tht_{0}(x)$.
Moreover let $(y_{0},\ldots,y_{\ell-1})$ be a list of all free variables occurring in $\pi$.
$Tm(\pi)$ denotes the set of all terms $t$, which is either the (induction) term $t$ in (\ref{eq:VJrule}),
or the (witnessing) term $t$ in 
\[
\infer[(\exi)]{\Gam}{\Gam,A(t)}
\]

Let $c(\pi)\geq 2$ be the number defined as follows.
First let $e_{1}(x)=x$, $e_{c+1}(x)=x^{e_{c}(x)}$.
Then $c=c(\pi)$ denotes a number for which the following holds for $e_{c}(x)$:
\beqn\label{eq:DHCon0}
\max\{t(y_{0},\ldots,y_{\ell-1}),\la y_{0},\ldots,y_{\ell-1}\ra : t\in Tm(\pi)\}\leq 
e_{c}(\max\{y_{0},\ldots,y_{\ell-1}\})
\eeqn
with a code $\la y_{0},\ldots,y_{\ell-1}\ra$ of the sequence $(y_{0},\ldots,y_{\ell-1})$.

By invoking the principle DH, let $K$ be a positive integer such that
\[
K\to_{\Del}\left(k+c+n-2\right)^{1+n}_{2^{m}}
\]
where $n=q(\pi)$, 
$m=m(\pi)+6$, $c=c(\pi)$,
$k=\max\{n, 3, 2 d(\pi)+5-c\}$
with
$K=\{0,1,\ldots,K-1\}$.

Let $Fml(\pi)=\{\vphi_{j}: j<m(\pi)\}$ be the set of all formulas occurring in $\pi$
other than $\Del_{0}$-formulas.
For formulas
$\vphi_{j}\equiv\vphi_{j}(y_{0},\ldots,y_{\ell-1})$,
let
$\vphi_{j}^{\prime}(y)\equiv\vphi_{j}((y)_{0},\ldots,(y)_{\ell-1})$.
Also for lists $X=(x_{1},\ldots,x_{n})$ of variables,
let
$\left(\vphi_{j}^{\prime}\right)^{(X)}(x_{0}):\equiv 
Q_{1}z_{1}<x_{1}Q_{2}z_{2}<x_{2}\cdots Q_{n}z_{n}<x_{n}
\tht((x_{0})_{0},\ldots,(x_{0})_{\ell-1})$.
Then $\Gam_{0}=\Gam_{0}(x_{0},x_{1},\ldots,x_{n})$ denotes the set of formulas
$\{\left(\vphi_{j}^{\prime}\right)^{(X)}(x_{0}): j<m(\pi)\}$ augmented with the six formulas
$\{c<x_{1}, x_{2}=x_{1}+x_{0}+1,x_{1}(x_{0}+1)<x_{2}, x_{1}^{x_{0}+1}<x_{2},
x_{1}^{x_{1}}<x_{2}, e_{c}(x_{1})<x_{2}\}$.

Let
$D=\{\alp_{-1}<\alp_{0}<\alp_{1}<\cdots<\alp_{k}<\alp_{k+1}<\cdots<\alp_{k+c+n-4}\}$
be a diagonal indiscernible subset of $K$ with respect to $\Gam_{0}$, 
cf.\,Proposition \ref{th:diagonalindiscernible}.

\bprp\label{eq:DHCon2}
Let $0\leq i<k$ and $\bet_{0},\ldots,\bet_{\ell-1}<\alp_{i}$.
Then
\\
$\max\{t(\bet_{0},\ldots,\bet_{\ell-1}), \la \bet_{0},\ldots,\bet_{\ell-1}\ra :t\in Tm(\pi) \}<\alp_{i+1}$.
\eprp
\bprf
First we see that $c<\alp_{1}$ from the fact that $D$ is indiscernible for $c<x_{1}$.
Second we see that
$0\leq i<j<p\leq k+3  \Rarw  \alp_{j}+\alp_{i}<\alp_{p}$,
$0\leq i<j<p\leq k+2 \Rarw  \alp_{j}\alp_{i}<\alp_{p}$, and 
$0\leq i<j<p\leq k+1 \Rarw  \alp_{j}^{\alp_{i}}<\alp_{p}$
from the indiscerniblity for
$x_{2}=x_{1}+x_{0}+1, x_{1}(x_{0}+1)<x_{2}, x_{1}^{x_{0}+1}<x_{2}$.
Third 
$0<i<j\leq k\Rarw \alp_{i}^{\alp_{i}}<\alp_{j}$ follows from the indiscernibility for
$x_{1}^{x_{1}}<x_{2}$.
Finally by the third and the indiscernibility for $e_{c}(x_{1})<x_{2}$ we obtain
$e_{c}(\alp_{i})<\alp_{i+1}$.
(\ref{eq:DHCon0}) yields the proposition.
\eprf
\\

For formulas $\sig\equiv\left(Q_{m}x_{m}\cdots Q_{n}x_{n}\tht\right)\,(1\leq m\leq n+1)$
and 
$I=\{\alp_{1}<\cdots<\alp_{n}\}\in[\Natural]^{n}$,
let $\sig^{(I)}$ be the $\Del_{0}$-formula
\[
\sig^{(I)}:\equiv\left(Q_{m}z_{m}<\alp_{1}\cdots Q_{n}z_{n}<\alp_{n-m+1}\tht\right)
\]
We write $\sig^{(I)}$ for $\Natural\models\sig^{(I)}$.

For any formula $\vphi$ occurring in the derivation $\pi$,
the following holds by the indiscernibiity for $\left(\vphi_{j}^{\prime}\right)^{(X)}(x_{0})$. 
{\small
\beqn\label{eq:DHCon3}
\bet_{0},\ldots,\bet_{\ell-1}< \alp_{i}<\alp_{i+1}<I,J\in[D]^{n}
\Rarw
\left[
\vphi^{(I)}(\bet_{0},\ldots,\bet_{\ell-1}) \Lrarw \vphi^{(J)}(\bet_{0},\ldots,\bet_{\ell-1}) 
\right]
\eeqn
}

To show Theorem \ref{th:DHCon}, it suffices to show
the following Theorem \ref{lem:DHCon}.
Let $\Gam$ be a sequent in the derivation $\pi$.
For
$J=\{\alp_{j},\alp_{j+1},\ldots, \alp_{j+n-1}\}\,(0<j\leq k)$,
let
$\vphi^{(j)}:\equiv\vphi^{(J)}$.

\bth\label{lem:DHCon}
Assume ${\sf PA}\vdash^{d}\Gam$, and for sequences
$(i_{0},\ldots,i_{\ell-1})\in \Natural^{\ell}$, let
$i=\max\{i_{j}: j<\ell\}\leq 2(d(\pi)-d)$.
Then
$\bigvee\{\vphi^{(i+2)}(\bet_{0},\ldots,\bet_{\ell-1}) : \vphi\in\Gam\}$ holds for
$\bigwedge_{j<\ell}(\bet_{j}<\alp_{i_{j}})$.
\end{theorem}

\bcor\label{cor:DHConint}
Assume ${\sf PA}\vdash Q_{1}z_{1}Q_{2}z_{2}\cdots Q_{n}z_{n}\tht$ with $\tht\in\Del_{0}$
and a sentence $Q_{1}z_{1}Q_{2}z_{2}\cdots Q_{n}z_{n}\tht$ for $Q_{i}\in\{\exi,\fal\}$.
Then there are natural numbers $\alp_{1}<\alp_{2}<\cdots<\alp_{n}$ such that
$\Natural\models Q_{1}z_{1}<\alp_{1}Q_{2}z_{2}<\alp_{2}\cdots Q_{n}z_{n}<\alp_{n}\tht$.
\ecor
{\bf Proof} of Theorem \ref{lem:DHCon} by induction on $d$.
Let us write $\vdash^{d}\Gam$ for ${\sf PA}\vdash^{d}\Gam$.

First consider the case for $(\exi)$.
\[
\infer[(\exi)]{\vdash^{d}\Gam}{\vdash^{d-1}\Gam,A(t)}
\]
Suppose $\lnot\bigvee\{\vphi^{(i+2)}(\bet_{0},\ldots,\bet_{\ell-1}) : \vphi\in\Gam\}$.
Then $\lnot\bigvee\{\vphi^{(i+3)}(\bet_{0},\ldots,\bet_{\ell-1}) : \vphi\in\Gam\}$ by $\vdash^{d}\Gam$.
IH yields for
$t_{1}\equiv t(\bet_{0},\ldots,\bet_{\ell-1})$,
$A^{(i+3)}(t_{1})$.
For the term $t\in Tm(\pi)$ we obtain $t_{1}<\alp_{i+1}<\alp_{i+2}$ by 
Proposition \ref{eq:DHCon2}.
Hence
$(\exi x A)^{(i+2)}\equiv(\exi x<\alp_{i+2}A^{(i+3)})$ follows.

Next consider the case for $(\fal)$.
\[
\infer[(\fal)]{\vdash^{d}\Gam}{\vdash^{d-1}\Gam,A(y)}
\]
By IH we can assume
$\bet<\alp_{i+2} \Rarw A^{(i+4)}(\bet)$.
Hence $\fal x<\alp_{i+2} A^{(i+4)}(x)$.
(\ref{eq:DHCon3}) yields $\fal x<\alp_{i+2} A^{(i+3)}(x)$, i.e.,
$(\fal x A)^{(i+2)}$.

Finally consider the case for the complete induction.
\[
\infer{\vdash^{d}\Gam}
{\vdash^{d-1}\Gam,A(0) 
& 
\vdash^{d-1}\Gam,\lnot A(a), A(a+1) 
& 
\vdash^{d-1}\Gam,\lnot A(t)}
\]
By IH we can assume $A^{(i+3)}(0)$.
For terms $t\in Tm(\pi)$ and $t_{1}\equiv t(\bet_{0},\ldots,\bet_{\ell-1})$,
we obtain $t_{1}<\alp_{i+1}$ by Proposition \ref{eq:DHCon2}.
IH yields
$\fal\bet<\alp_{i+1}\left[A^{(i+3)}(\bet) \to A^{(i+3)}(\bet+1)\right]$, and
$A^{(i+3)}(t_{1})$.
On the other hand we have $\lnot A^{(i+3)}(t_{1})$ by IH.
\eprf
\\

Assume ${\sf PA}\vdash^{d(\pi)}\exi x\,\tht_{0}\,(\tht_{0}\in\Del_{0})$.
Theorem \ref{lem:DHCon} yields 
$\left(\exi x\,\tht_{0}\right)^{(2)}\equiv\exi x<\alp_{2}\,\tht_{0}$.
Theorem \ref{th:DHCon} is shown.
\\

For positive integers $n,m,k$, let
\[
D(n,m,k):=\min\{K\geq\max\{n+2,m\}: K\to_{\Del}(k)_{m}^{1+n}\}
\]

\bcor\label{cor:DHCon}
\benu
\item\label{cor:DHCon.1}
Let $f(x)$ be a provably computable function in {\sf PA}.
Then there exists an $n_{0}$ such that
$\fal x\in\Natural\left(f(x)<D_{n_{0}}(x)\right)$ for
$x\mapsto D_{n_{0}}(x)=D(n_{0},n_{0},n_{0}+x)$.

\item\label{cor:DHCon.2}
$n\mapsto D(n)=D(n,n,2n)$ dominates every provably computable function in {\sf PA}.
\eenu
\ecor
\bprf
\ref{cor:DHCon}.\ref{cor:DHCon.1}.
Let $\pi=\pi(x)$ be a derivation of $\exi y\,\tht(x,y)$ in {\sf PA}, and
$n_{0}=\max\{2^{m(\pi)+6},2q(\pi)+c(\pi)-1, 3+c(\pi)-1, q(\pi)+2d(\pi)+5\}$.
For a natural number $a$, let $\Gam_{a}$ be a set of formulas obtained from
$\Gam_{0}$ by replacing the formula $c<x_{1}$ by $\max\{c,a\}<x_{1}$.
Let $K_{a}=D_{n_{0}}(a)$, and
$D_{a}=\{\alp_{-1}<\alp_{0}<\alp_{1}<\cdots\}$ be a diagonal homogeneous subset of
$K_{a}$ for $\Gam_{a}$, and in size 
$\#(D_{a})\geq n_{0}+a$.
We see easily $a<\alp_{1}$.
$\pi(a)$ denotes the derivation of $\exi y\,\tht(a,y)$ in {\sf PA}
obtained from $\pi(x)$ by substituting the numeral $a$ for the variable $x$.

Theorem \ref{lem:DHCon} yields $\exi y<\alp_{3}\,\tht(a,y)$,
and $\exi y<K_{a}\,\tht(a,y)$ by $\alp_{3}<K_{a}$,
\\
\ref{cor:DHCon}.\ref{cor:DHCon.2}.
This follows from Corollary \ref{cor:DHCon}.\ref{cor:DHCon.1} and the fact
$D_{n_{0}}(a)\leq D(a)$ for $a\geq n_{0}$.
\eprf

\end{document}